\documentclass{amsart}
\usepackage{amssymb}  
\usepackage{amsthm}
\usepackage{amsmath,latexsym,amssymb} 
\usepackage[all]{xy}
\usepackage{graphics}
\usepackage{graphicx}
\usepackage{hyperref}
\newcommand{\rp}[1]{\mathbb{RP}^{#1}} 
\newcommand{\map}{\mbox{map}}
\newcommand{\omcg}{\Gamma}
\newcommand{\nmcg}{\mathcal{N}}

\newcommand{\us}[1]{\underset{#1}{\times}}
\newcommand{\res}[1]{\textit{#1}}
\newcommand{\diff}{\mbox{Diff}}
\newcommand{\topp}{\mbox{Top}}
\newcommand{\br}{\mathbb{R}}

\newcommand{\barra}{\;|\;}
\newcommand{\baca}{\backslash}
\newcommand{\bas}{\backslash}

\newcommand{\bz}{\mathbb{Z}}
\newcommand{\ff}{\mathbb{F}}
\newcommand{\ra}{\rightarrow}
\newcommand{\lra}{\longrightarrow}

\newcommand{\stac}[1]{\stackrel{#1}{\longrightarrow}}
\newcommand{\til}[1]{\widetilde{#1}}
\newcommand{\wt}{\widetilde}
\newcommand{\klein}{\mathbb{K}}

\newtheorem{theorem}{Theorem}[section]
\newtheorem{lemma}[theorem]{Lemma}

\theoremstyle{definition}

\theoremstyle{remark}

\numberwithin{equation}{section}

\begin{document}

\title{Mapping class group and \\function spaces: a survey}


\author{Fred R. Cohen}
\address{{\it Department of Mathematics}\\
University of Rochester\\
Rochester\\
NY 14627}
\curraddr{}
\email{cohf@math.rochester.edu}
\thanks{}

\author{Miguel A. Maldonado}
\address{{\it Unidad Acad\'emica de Matem\'aticas}\\
Universidad Aut\'onoma de Zacatecas\\
Zacatecas\\
M\'exico  98000}
\curraddr{}
\email{mimaldonadoa@gmail.com}
\thanks{}

\subjclass[2010]{Primary 55P20, 20F36, 37E30}

\keywords{Mapping class groups, function spaces.}

\date{October 22, 2014}

\dedicatory{This article is dedicated to Samuel Gitler Hammer\\ who brought much joy and interest to mathematics.}

\begin{abstract}
This paper is a survey of the relationship between labelled configuration spaces, mapping class groups with marked points and function spaces. In particular, we collect calculations of the cohomology groups for the mapping class groups of low-genus orientable and non-orientable surfaces.
\end{abstract}

\maketitle

\tableofcontents

\section{Introduction}\label{sec1}
The \res{mapping class group} $\omcg(S)$ of an orientable surface (closed, connected, orientable compact smooth 2-manifold) $S$ is the group of isotopy classes in $\topp^+(S)$, the group of orientation-preserving homeomorphisms of $S$ or, in other words, $\omcg(S)=\pi_0\topp^+(S)$ \cite{Birman, FM}. Since we can define the mapping class group within other contexts (PL category, smooth category) obtaining isomorphic groups, in various parts of this paper we sometimes use the diffeomorphism group $\diff^+(S)$ instead of $\topp^+(S)$ \cite{Boldsen09}. 

\medskip
The mapping class group of a surface is closely related to many parts of mathematics. Here, we present the relationship of this group to the homotopy theory of configuration spaces. In particular, we are interested in how configuration spaces of surfaces are used in (co)homology calculations for mapping class groups by considering these spaces and some generalizations in the construction of Eilenberg-Mac Lane spaces for $\omcg(S)$ as well as some variants of it. This relation allows us to consider the unstable range of the homology and cohomology of mapping class groups by adopting a homotopical point of view. In particular we use the function space $\map(X,Y)$, and pointed function space $\map_*(X,Y)$ for obtaining ``global information'' for the (co)homology of mapping class groups with all finite marked points.

\medskip

The authors thank Ulrike Tillmann for her useful suggestions. 

\section{Mapping class groups with marked points}\label{sec2}
There are variants of the mapping class group $\omcg(S)$ which are  determined by restricting to certain substructures on $S$ and/or by restrictions on the behaviour of homeomorphisms on that structure. The main case in the present survey is the \res{mapping class group with marked points} $\omcg^k(S)$ defined as the group of isotopy classes in $\topp^+(S;k)$, the group of homeomorphisms which leave a set of $k$ distinct points in $S$ invariant, that is $\omcg^k(S)=\pi_0\topp^+(S;k)$. We define the \res{pure} mapping class group with $k$ marked points $P\omcg^k(S)$ as the kernel of the natural homomorphism $\omcg^k(S)\ra \Sigma_k$, obtained by the action of $\omcg^k(S)$ on punctures. For the orientable surface $S_g$ of genus $g$ we will use the following notation  $\omcg_g=\omcg(S_g)$ and $\omcg_g^k=\omcg^k(S_g)$. Also, when one of the parameters $g,k$ is zero it will be omitted from the notation.

\medskip
One reason for considering marked points is that information about low genus surfaces gives information about surfaces of higher genus via the theory of branched covers. One example is the central extension $1\ra \bz/2\bz\ra \omcg_2\ra \omcg^6\ra 1$ of \cite{BergauMennicke60,BirmanHilden71}  used in \cite{BensonCohen91} to work out the homology of $\omcg_2$ from that of $\omcg^6$. Also, for genus $g\geq2$ surfaces, the well-known Birman short exact sequence $1\ra \pi_1F_k(S_g)/\Sigma_k \ra \omcg_g^k \ra \omcg_g\ra 1$
is used in inductive processes; see for example \cite{Lu02} where this sequence is used to inherit cohomological properties from $\omcg_g$ to mapping class groups with marked points.

\medskip
The homology and cohomology of a discrete group $G$ can be studied from a homotopical viewpoint by considering Eilenberg-Mac Lane spaces $K(G,1)$. For the braid group of a surface different from the sphere $S^2$ and the projective plane $\rp{2}$ the configuration spaces $F_k(M)$ and $F_k(M)/\Sigma_k$, are $K(G,1)$'s, where the fundamental groups are denoted $P_k(S)$ and $B_k(S)$, respectively. In the case the surface $S$ is $\br^2$, then the groups $P_k(\br^2)$ and $B_k(\br^2)$ are Artin's pure braid group with $k$ strands, respectively Artin's braid group with $k$ strands. For the cases of $S^2$ and $\rp{2}$ certain Borel constructions are needed. This setting can be considered as a starting point for the results exposed here. For completeness we recall this construction in what follows.

\medskip
Consider the group $\topp^+(S)$ acting transitively on the configuration space $F_k(S)/\Sigma_k$. Note that for a basepoint  $\hat{x}\in F_k(S)/\Sigma_k$, the isotropy subgroup is precisely $\topp^+(S;k)$. From this setting, there is an induced homeomorphism $\topp^+(S)/\topp^+(S;k)\cong F_k(S)/\Sigma_k$.  Under these conditions, the homotopy orbit space
\begin{equation}\label{borel}
E\topp^+(S)\us{\topp^+(S)}F_k(S)/\Sigma_k
\end{equation}
is homotopy equivalent to $E\topp^+(S)/\topp^+(S;k)$, a model for the classifying space of $\topp^+(S;k)$. Note that the fundamental group of this space is isomorphic to $\pi_1 B\topp^+(S;k)\cong \omcg^k(S)$, the mapping class group of $S$ with marked points. 

\medskip
Classical theorems concerning the homotopy type of the group \linebreak $\diff^+(S)$ (or $\topp^+(S)$, see for example \cite{EarleEells69,Gramain73}) give a better description of the space (\ref{borel}) above allowing applications to compute cohomology for several cases. These cases will be recalled in the next sections relating the spaces involved with function spaces. For example, if $k \geq 3$, then the construction $ESO(3)\times_{SO(3)}F(S^2,k)/\Sigma_k$
is a $K(\Gamma^k,1)$ (\cite{BensonCohen91,Cohen93, BodigheimerCohenPeim01}).

\section{On the labelled configuration space}\label{sec3}
In this section we will recall the main properties of the labelled configuration space $C(M;X)$. All the theory presented here is an exposition of the paper \cite{BodigheimerCohenTaylor89} as well as \cite{Cohen96}.

\medskip
For a smooth $m$-manifold $M$ we consider the configuration space with \res{labels in $X$} as the space $C(M;X)$ of finite configurations of points in $M$ parametrized by $X$, where $X$ is pointed space with base-point denoted by $*$.
\[
C(M;X)=\left.\left(\coprod_{k\geq 1}F_k(M)\us{\Sigma_k} X^k \right)\right/\approx,
\]
under the relation $\approx$ given by
\[
[m_1,\ldots,m_k;x_1,\ldots,x_k]\approx [m_1,\ldots,m_{k-1};x_1,\ldots,x_{k-1}], \;\mbox{if}\;\;x_k=*.
\]
Equivalently, the space $C(M;X)$ consist of equivalence classes of pairs $[S,f]$, where $S\subset M$ is a finite set and $f:S\ra X$. The equivalence relation is generated by
$$
[S\bas \{k\},f|_{S\bas \{k\}}]\approx [S,f], \;\;\; \mbox{for}\;f(k)=*
$$
Note that for $X=S^0$ the space $C(M;X)$ is the union of all finite unordered configurations of $M$.

\medskip
There is a description of $C(M;X)$ as a space of sections that will be used throughout the present exposition. In what follows we will recall this relation. Consider the disc bundle $D(M)\subset T(M)$ in the tangent bundle $T(M)$ of $M$ and define $\dot{T}(M)$ by the identification of two vectors that lie in the same fiber and outside $D(M)$, this gives a $S^m=D^m/\partial D^m$-bundle $\dot{T}(M)\ra M$. Denote by $\dot{T}(M;X)$ the fibrewise smash product of the bundle $\dot{T}(M)\ra M$ with $X$ to obtain a bundle
$$
\dot{T}(M;X)\lra M
$$
where $\dot{T}(M;X)$ is the total space over $M$ with fibre $\Sigma^m X$. 

\medskip
Next, consider $\Gamma(M;X)$ the space of sections of $\dot{T}(M;X)$. In the special case of $X=S^0$, D. McDuff (\cite{McDuff75}) gives a proof that there is a map $C(M;S^0)\ra \Gamma(M;S^0)$ which induces a homology isomorphism through a range increasing with $k$ on the level of a map $F_k(M)/{\Sigma_k} \to \Gamma_k(M;S^0), \  \mbox{the path-component of degree $k$ maps}$.  A further statement in case $X$ is path-connected is given in \cite{Bodigheimer85} and \cite{Francis14}: There is a homotopy equivalence
$$
\gamma:C(M;X)\lra \Gamma(M;X)
$$
for path-connected $X$. 

\medskip
In particular, if $M$ is parallelizable then $C(M;X)\simeq \map(M,\Sigma^mX)$. Moreover, one has $\Gamma(M;\Sigma^kX)\simeq \map(M,\Sigma^{m+k}X)$, where the tangent bundle of $M$ plus a trivial $k$-plane bundle over $M$ is a trivial bundle. Thus for example if $M=S^n$, then the tangent bundle of $S^n$ plus a trivial line bundle gives a trivial $(n+1)$-plane bundle.

\medskip
The case when the manifold is a sphere will be of importance for relating function spaces with mapping class groups so we recall it here:

\begin{lemma}[\cite{BodigheimerCohenTaylor89}]\label{maps}
For a connected CW-complex $X$ there is a homotopy equivalence
$$
C(S^m;\Sigma X)\simeq \mbox{\em map}(S^m,\Sigma^{m+1}X)=\Lambda^m\Sigma^{m+1}X.
$$
\end{lemma}

The homotopy equivalence $\gamma$ above commutes with the map $\phi^*$ induced by an isometry $\phi:M\ra M$ which acts on $C(M;X)$ and $\Gamma(M;X)$ in a natural way inducing the following commutative diagram
$$
\xymatrix{C(M;X)\ar[d]^{\phi^*}\ar[r]^\gamma&\Gamma(M;X)\ar[d]^{\phi^*}\\C(M;X)\ar[r]^\gamma&\Gamma(M;X).}
$$
Moreover, assume that $M$ is a Riemannian manifold with a subgroup $G$ of the isometry group of $M$. Then there is a $G$-equivariant homotopy equivalence (\cite{Bodigheimer.Madsen,Manthorpe.Tillmann}).
\begin{equation}\label{borelsection}
EG\times_GC(M;X)\lra EG\times_G\Gamma(M;X)
\end{equation}
which will be used later for relating the homotopy orbit space of the previous section with explicit calculations for mapping class groups with marked points.

\medskip
The space $C(M;X)$ has a natural filtration $C_k(M;X)$ given by lengths of configurations 
and for  $D_k(M;X)=C_k(M;X)/C_{k-1}(M;X)$ the space $C(M;X)$ stably splits as the wedge $\bigvee_{k\geq1} D_k(M;X)$ (\cite{CohenMayTaylor84}, \cite{Cohen83}) and so there is an homology isomorphism 
$$
\wt{H}_*(C(M;X);\ff)\cong \bigoplus_{k=1}^\infty \wt{H}_*(D_k(M;X);\ff)
$$
for any field $\ff$ of characteristic zero.

\medskip
Daniel Kahn first gave a stable decomposition for $\Omega^{\infty}\Sigma^{\infty}(X)$ \cite{DanielKahn}. V. Snaith first proved that $\Omega^n \Sigma^n(X)$ stably splits \cite{Snaith}. Subsequently, it was shown that a similar decomposition
applies to $C(M,X)$ in case $X$ has a non-degenerate base-point, and $M$ is any "reasonable" space.
\cite{CohenMayTaylor84,Cohen83} as observed earlier in \cite{CohenTaylor76}.

\subsection{Some notes on the homology of $F_k(M)/\Sigma_k$}\label{sub3.1}
In the case of configurations of points parametrized by a sphere $S^n$ the space $D_k(M;S^n)$ is equal to $F_k(M)\times_{\Sigma_k} S^{nk}/F_k(M)\times_{\Sigma_k} *$. The action of $\Sigma_k$ on $H_{nk}(S^{nk})=\bz$, with the usual sign conventions, is trivial if $n$ is even and is the sign action if $n$ is odd and this information forces the homology of $D_k(M;S^n)$ to give the homology of $F_k(M)/\Sigma_k$ if $n$ is even and if $n$ is odd it gives the homology of $F_k(M)/\Sigma_k$ with coefficients in the sign representation. More information will be given in Section 4 below and for the rest of the section we will assume that $n$ is even.

\medskip
The space $D_k(M;S^n)$ is homeomorphic to the Thom space $T(\eta)$ of a vector bundle over the configuration space $F_k(M)/\Sigma_k$ (\cite{Wang02}) giving an isomorphism for a fixed $k$
\begin{eqnarray*}
H_*(F_k(M)/\Sigma_k;\ff)&\cong& H_{*+kn}(D_k(M;S^n);\ff)\\
&\cong& H_{*+kn}(C(M;S^n);\ff)
\end{eqnarray*}
for $0<j<mk$ where $m$ is the dimension of $M$ with sufficiently large $n$.

From this result it is possible to describe the homology of $F_k(M)/\Sigma_k$ additively since the homology of $C(M;S^n)$ is given as the graded vector space

\vspace{-4mm}
\begin{equation}\label{graded}
\bigotimes_{q=0}^m \bigotimes_{\;\;\beta_q} H_*(\Omega^{m-q}S^{m+n};\ff)
\end{equation}
where $\beta_q$ is the $q$-th Betti number of $M$ (that is, the dimension of $H_q(M;\ff)$) and $\ff$ is a field of characteristic zero with the restriction of $m+n$ to be odd if computations are taken over other field than $\ff_2$ (\cite{BodigheimerCohenTaylor89}). Every factor on the product (\ref{graded}) above is an algebra with weights associated to its generators and the reduced homology of $D_k(M;S^n)$ is the vector subspace generated by all elements of weight $k$. Finally, the additive structure of the homology of unordered configurations can be found by counting generators of weight $k$ (See \cite{BodigheimerCohenTaylor89}). 

\medskip
When $M$ is a surface the mod-$2$ homology of the labelled configuration space $C(M;S^n)$ is given by the tensor product $
H_*(\Omega^2S^{n+2})^{\otimes \beta_0}\otimes H_*(\Omega S^{n+2})^{\otimes \beta_1}\otimes H_*(S^{n+2})^{\otimes \beta_2}$. We will use the term for $\Omega^2S^{n+2}$ to express homology of unordered configurations in terms of the homology of classical braid groups. 

\medskip
First consider the product (\ref{graded}) above when $M$ is a surface with Betti numbers $\beta_0,\beta_1,\beta_2$ given by $\beta_0=1=\beta_2$, and $\beta_1=2g$ for mod-$2$ homology. Then
\[
H_*(C(M;S^n);\ff_2)\cong H_*(\Omega^2S^{n+2};\ff_2)\otimes (\ff_2[x])^{\otimes 2g}\otimes (\ff_2[u]/u^2).
\]
An element on this product has the form $y\otimes x^l\otimes u^{\epsilon}$, where $y$ runs over an additive basis for $H_*(\Omega^2S^{n+2};\ff_2)$ and $x, u$ are the fundamental classes on degrees $n+1$ and $n+2$, respectively. The weight function is given as
$$
\omega(u)=1,\;\; \omega(x^i)=i,\;\;\;\mbox{for}\; i=1,2,\ldots
$$
The generator $y$ relates this setting with the braid groups $B_k$ by recalling that $C(\br^2;S^n)$ is homotopy equivalent to $\Omega^2\Sigma^2S^n=\Omega^2S^{n+2}$ from which one gets 
\begin{eqnarray*}
\til{H}_q(\Omega^2 S^{n+2};\ff_2)&\cong& \bigoplus_k \til{H}_{q-kn} (F_k(\br^2)/\Sigma_k;\ff_2)\\
&\cong& \bigoplus_k \til{H}_{q-kn} (B_k;\ff_2)
\end{eqnarray*}

\medskip
Similar methods were developed in \cite{CohenLadaMay} where the equivariant homology of configurations $F_k(\br^n)$ with coefficients in the following modules were first worked out:

\begin{enumerate}
\item The trivial representation,

\item the sign representation, 

\item representations afforded by left cosets of Young subgroups of the symmetric groups, and 

\item the homology of the space $F_k(\br^n)\times_{\Sigma_k}X^k$ where $X$ is any compactly generated weak Hausdorff space with non-degenerate base-point (\cite{CohenLadaMay} page $226$, Corollary $3.3$).
\end{enumerate} 

These results were given in terms of analogues of graded Lie algebras with certain operations. The methods use $\br^n$ in basic ways.

\subsubsection{Configurations on a sphere}
Let us move a step backwards by letting $M=S^r$ and consider the isomorphism
\begin{equation}\label{sphere}
H_*(C(S^r;S^n);\ff_2)\cong H_*(\Omega^rS^{r+n};\ff_2)\otimes \ff_2[u]/u^2.
\end{equation}
where the generators in this product have the form $y\otimes u^\epsilon$, for $\epsilon=0,1$.

\medskip
As in the last page, we use that $C(\br^r;S^n)$ is homotopy equivalent to $\Omega^r\Sigma^r S^n\simeq \Omega^rS^{r+n}$ to get
$$
\til{H}_q(\Omega^r S^{r+n};\ff_2)\cong \bigoplus_k \til{H}_{q-kn} (F_k(\br^r)/\Sigma_k;\ff_2)
$$
Since the homology of $F_k(S^r)/\Sigma_k$ is obtained from that of $C(\br^r;S^n)$ by considering monomials of weight $k$ we have the generator $y$ must have degree $|y|=q+kn-\epsilon(r+n)$ and weight $\omega(y)=k-\epsilon$ according to which one has
\begin{eqnarray*}
H_{q+rn}(D_k(\br^r;S^n);\ff_2)&\cong& H_q (F_k(\br^r)/\Sigma_k;\ff_2)\\
H_{q+rn-(r+n)}(D_{k-1}(\br^r;S^n)&\cong& H_{q-r} (F_{k-1}(\br^r)/\Sigma_{k-1};\ff_2)
\end{eqnarray*}
From these expressions we finally get the isomorphism of vector spaces
$$
H_q(F_k(S^r)/\Sigma_k;\ff_2)\cong H_q(F_k(\br^r)/\Sigma_k;\ff_2) \bigoplus H_{q-r}(F_{k-1}(\br^r)/\Sigma_{k-1};\ff_2)
$$
which in the case of the 2-sphere gives
\begin{eqnarray*}
H_q(F_k(S^2)/\Sigma_k;\ff_2)\cong H_q(B_k;\ff_2) \bigoplus H_{q-2}(B_{k-1};\ff_2)
\end{eqnarray*}
as stated in \cite{BodigheimerCohenPeim01}.

\subsubsection{Configurations on non-orientable surfaces}
Let $N_g$ be the compact non-orientable surface of genus $g$, where genus means the number of projective planes in a connected sum decomposition; that is, $N_1$ is $\rp{2}$ and $N_2$ is the Klein bottle. The mod-$2$ homology of labelled configurations in $N_g$ is given as the tensor product
\[
H_*(C(N_g;S^n);\ff_2)\cong H_*(\Omega^2S^{n+2};\ff_2)\otimes (\ff_2[x])^{\otimes g}\otimes \ff_2[u]/u^2.
\]
where $x, u$ are the fundamental classes on degrees $n+1$ and $n+2$, respectively. In what follows we will describe the relation of this isomorphism with braid groups for the cases $g=1,2$. More details can be found in \cite{MaldonadoXicotencatl14}.

\medskip
For the projective plane $\rp{2}$ consider the tensor product 
\[
H_*(C(\rp{2};S^n))\cong H_*(\Omega^2S^{n+2}) \otimes \ff_2[x]\otimes \ff_2[u]/u^2,
\]
with basis monomials of the form $y\otimes x^r \otimes  u^\epsilon$, for $\epsilon =0,1$ and $r =0, 1, 2, \dots$. Degree and weight for elements in the form $x^r\otimes u^\epsilon$ are given by $\epsilon n + 2\epsilon + rn + r$ and $\epsilon+r$, respectively. Thus, fixing $q$ and  $k$, the monomial  $y\otimes x^r \otimes u^\epsilon$ represents a generator in $H_q  (F_k(\rp{2})/\Sigma_k;\ff_2)    \cong    H_{q + kn}(D_k(\rp{2}; S^n);\ff_2$ if and only if the degree and weight of $y$ are given by $|y|=q + (k-r-\epsilon)n - 2\epsilon -r$ and $\omega(y) =k  - (\epsilon + r)$. A careful analysis on the possibilities of the numbers involved lead to the following isomorphisms of underlying spaces
\begin{eqnarray*} 
H_{q+(k-r)n-r}  \left( D_{k-r} (\br^2;  S^n) ;\ff_2  \right)  
\cong H_{q-r}(B_{k-r};\ff_2)
\end{eqnarray*}
and
\begin{eqnarray*} 
H_{q+(k-r-1)n-r-2}  \left( D_{k-r-1} (\br^2;  S^n) ;\ff_2 \right)
\cong H_{q-r-2}(B_{k-r-1};\ff_2)
\end{eqnarray*}
Finally we get that $H_q  (F_k(\rp{2})/\Sigma_k;\ff_2)$ is isomorphic to
$$
\bigoplus_{r=0}^{\text{min}\{q,k\}} H_{q-r}(B_{k-r};\ff_2)   
\;  \oplus  \;    \bigoplus_{\ell=0}^{\text{min}\{q-2,k-1\}}  H_{q-r-2}(B_{k-r-1};\ff_2)
$$
as graded vector spaces.

\medskip
In the case of configurations on the Klein bottle $\klein$ one considers the isomorphism
$$
H_*(C (\klein; S^n)) \cong   H_*(\Omega^2 S^{n+2} )\otimes \ff_2[x_1, x_2] \otimes \ff_2 [u]/u^2, 
$$
having basis all elements of the form 
$u^\epsilon \otimes x_1^r \otimes  x_2^s \otimes y$, for $\epsilon =0,1$ and $r,s =0, 1, 2, \dots$. Notice that monomials of the form $u^\epsilon\otimes x_1^r \otimes x_2^s$ have degree and weight given by
\begin{eqnarray*}
  |u^\epsilon\otimes x_1^r\otimes x_2^s| &=& \epsilon n + 2\epsilon + (r+s)(n+1),\\
   \omega(u^\epsilon\otimes x_1^r \otimes x_2^s) &=&  \epsilon + r + s
\end{eqnarray*}
  
\noindent
As before, fix $q$ and  $k$ and note the degree and weight of $y$ are given by
\begin{align*}
|y|   &  =  q + (k-r-s-\epsilon)n - 2\epsilon - r-s\\
\omega(y) &= k  - (\epsilon + r +s) =  k -\epsilon - r-s
\end{align*}
leading to the isomorphisms
\begin{align*} 
H_{q + (k-r-s)n - r-s} \left(  D_{k-r-s} (\br^2;  S^n) ;\ff_2 \right)
\cong H_{q-r-s}(B_{k-r-s};\ff_2)
\end{align*}
$$
H_{q + (k-r-s-1)n -r-s-2}  ( D_{k-r-s-1} (\br^2;  S^n) ;\ff_2 )   
\cong
$$ 
$$
\cong
H_{q-r-s-2}(B_{k-r-s-1};\ff_2)
$$	

Finally, from these expressions, $H_q  (F_k(\klein)/\Sigma_k;\ff_2)$ is isomorphic as graded vector space, to the direct sum
$$
\bigoplus_{r,s} H_{q-r-s}(B_{k-r-s};\ff_2)   \;\;  \oplus  \;\;    \bigoplus_{r,s} H_{q-r-s-2}(B_{k-r-s-1};\ff_2),
$$
where the indices $r,s$ run over the non-negative integers such  $r+s  \leq  \text{min}\{q,k\}$ in the first direct sum and $r+s  \leq  \text{min}\{q-2,k-1\}$ in the second direct sum.

\section{A bundle-generalization of $C(M;X)$}\label{sec4}
We collect here some results on a bundle-generalization of the labelled configuration space $C(M;X)$ taken from \cite{BodigheimerCohenPeim01}. Let $\pi:E\ra B$ be a fibre bundle with fibre $Y$ and define
$$
E(\pi,k)=\{(x_1,\ldots,x_k)\in E^k \barra x_i\neq x_j, \pi(x_i)=\pi(x_j), i\neq j\}.
$$ 
The space $E(\pi,k)$ and the unordered version $E(\pi,k)/\Sigma_k$ are total spaces of fiber bundles with fiber configuration spaces
$$
F_k(Y)\lra E(\pi,k)\stac{p} B,\qquad F_k(Y)/\Sigma_k\lra E(\pi,k)/\Sigma_k\stac{p} B
$$
where $p$ is the projection on the first coordinate. These fiber-analogues of configuration spaces have appeared in \cite{Cohen93}. As before, let $X$ be a CW-complex with basepoint $*$ and define
\vspace{-1.5mm}
\[
E(\pi;X)=\left(\coprod_{k\geq 0}E(\pi,k)\us{\Sigma_k} X^k \right)/\approx,
\]
where $\approx$ is the relation given for $C(M;X)$ by requiring $\pi(S)$ to be a single point. In fact, note that when $B$ is a point, then $E(\pi,X)=C(Y,X)$ and so $E(\pi;X)$ are bundle-generalizations of labelled configuration spaces.

\medskip
There is also a description of $E(\pi;X)$ in terms of pairs $[S,f]$ requiring that $\pi(S)$ is a single point. So there is a filtration $E_k(\pi;X)$ by the cardinality of $S$. As before, if we denote by $D_k(\pi;X)$ the filtration quotient $E_k(\pi;X)/E_{k-1}(\pi;X)$ then $E(\pi;X)$ is stably equivalent to the wedge $\bigvee_{k\geq 1} D_k(\pi;X)$ and thus there is a homology isomorphism (\cite{Cohen83})
$$
\wt{H}_*(E(\pi;X);\ff)\cong \bigoplus_k \wt{H}_* (D_k(\pi;X);\ff)
$$

Let $M$ be a $G$-space and consider the Borel fibration $\lambda: EG\times_{G}M\ra BG$ with fibre $M$. There are homeomorphisms 
$$
E(\lambda,k)\cong EG\times_G F_k(M),\qquad E(\lambda,k)/\Sigma_k\cong EG\times_G F_k(M)/\Sigma_k
$$
and so $E(\pi,k),E(\pi,k)/\Sigma_k$ are bundle-generalizations of Borel constructions and they are used for (co)homology calculations for mapping class groups. It is also convenient to consider the Borel construction for labelled configuration spaces and so we have a homeomorphism
$$
E(\lambda;X)\cong EG\times_G C(M;X)
$$
where the $G$-action on $C(M;X)$ is induced from that on $M$. This homeomorphism is compatible with filtrations in the sense that there is an induced homeomorphism $E_k(\lambda;X)\cong EG\times_G C_k(M;X)$, where $C_k(M;X)$ is the $k$-th filtration of $C(M;X)$. Moreover, for the singular chain complex $S_*(X)$ we have a homology isomorphism $H_*(X;\ff)\cong S_*(X)$; in particular, there is an isomorphism
$$
S_*(E(\eta;k))\otimes_{\Sigma_k} H_*(X)^{\otimes k}\cong S_*(E(\eta;k)\times_{\Sigma_k} X^k)
$$

If $X^{(k)}$ denotes the $k$-fold smash product of $X$ then we get
$$
S_*(E(\eta;k))\otimes_{\Sigma_k} H_*(X)^{\otimes k}\cong S_*(\frac{E(\eta;k)\times_{\Sigma_k} X^{(k)}}{E(\eta;k)\times_{\Sigma_k} *})\cong S_*(D_k(\eta;X))
$$

Recall that for a graded vector space $V$ over the a field $\ff$ there exists a bouquet of spheres $S_V$ such that $V=\wt{H}_*(S_V;\ff)$. Let $X=S_V$ to get $S_*(E(\eta;k))\otimes_{\Sigma_k} V^{\otimes k}\cong S_* (D_k(\eta;S_V))$; that is
\begin{equation}\label{www}
H_*(E(\eta,k);V^{\otimes k})\cong \wt{H}_*(D_k(\eta;S_V);\ff)
\end{equation}
Note that the symmetric group $\Sigma_k$ acts on $V^{\otimes k}$ by permutation with the usual sign conventions. From this, it follows that if $V$ is a copy of $\ff$ concentrated in even degrees then $V^{\otimes k}$ is the trivial $\ff\Sigma_k$-module. On the other hand, if $V$ is concentrated in odd degrees, then $V^{\otimes k}$ is a copy of the sign representation.

\medskip
The isomorphism (\ref{www}) above will be used through the next sections for the bundle $BSO(n)\ra BSO(n+1)$ with fibre $S^n$, induced by the inclusion $SO(n)\subset SO(n+1)$.

\section{Cohomology of mapping class groups} \label{sec5}
There are many motivations for considering the (co)homology of mapping class groups but from a homotopical point of view the main motivation is that for $g>2$ it turns out that $B\topp^+(S_g)$ is the classifying space of $\omcg(S_g)$ and thus characteristic classes of $S_g$-bundles are classes in $H^*(\omcg(S_g);\bz)$ (\cite{Morita01}).

\medskip
In this section we collect the material concerning the spaces $C(M;X)$ and $E(\eta;X)$ related to cohomological calculations for mapping class groups of surfaces with marked points.

\subsection{The mapping class group of the sphere with marked points}
It is a classical result (\cite{FadellNeuwirth62}) that for $k\geq 3$, there is a homeomorphism $F_k(S^2)=SO(3)\times F_{k-3}(S^2\backslash Q_3)$ and thus $F_k(S^2)$ is never a $K(\pi,1)$. The Eilenberg-Mac Lane space $K(\pi,1)$ is given by the Borel construction (\ref{borel}) given in Section \ref{sec2}: for $k\geq 3$, this space is a $K(\omcg^k,1)$ space (\cite{Cohen93}). Considering the $SO(3)$-action on $S^2$ by rotations there is an induced diagonal action on $F_k(S^2)/\Sigma_k$. Moreover, since the inclusions $SO(3)\subset \diff^+(S^2)\subset \topp^+(S^2)$ are homotopy equivalences (\cite{Smale}) the Borel construction mentioned above is homotopy equivalent to 
$$
ESO(3)\times_{SO(3)}F_k(S^2)/\Sigma_k
$$
This space can be stated in terms of the combinatorial models given in the last section. For the bundle $\eta: BSO(2)\ra BSO(3)$ with fiber $S^2$ and $k\geq 3$ the space $E(\eta,k)/\Sigma_k$  is $K(\omcg^k,1)$ (\cite{BodigheimerCohenPeim01}); if $k=0,1,2$ the spaces are not $K(\pi,1)$ and are exceptional. Thus one has, for $k\geq 3$, an isomorphism $H^*(\omcg^k;R)\cong H^*(E(\eta,k)/\Sigma_k;R)$ for a trivial $\omcg^k$-module $R$. Moreover, from the isomorphism (\ref{www}) on the last section we have
\begin{eqnarray}
H_*(\omcg^k;\ff)\cong& \!\!\!\!\!\!\!\!\!\!\!\!\wt{H}_{*+2qk}(D_k(\eta;S^{2q});\ff)\\
H_*(\omcg^k;\ff(\pm 1))\cong& \wt{H}_{*+(2q+1)k}(D_k(\eta;S^{2q+1});\ff)
\end{eqnarray}

\medskip
These isomorphisms are considered in \cite{BodigheimerCohenPeim01} for giving the mod-$2$, mod-$p$ cohomology as also for coefficients in the sign representation. In what follows we recall the methods considered there. First consider the homotopy equivalence (\ref{borelsection}) for $M=S^n$ and $G=SO(n)$ to get
$$
ESO(n)\times_{SO(n)}C(S^n;\Sigma X)\simeq ESO(n)\times_{SO(n)}\Gamma(S^n;\Sigma X)
$$
The left side is simply the space $E(\eta;\Sigma X)$ with $BSO(n)\stac{\eta} BSO(n+1)$ and thus we get a homotopy equivalence
$$
E(\eta,\Sigma X)\simeq ESO(n)\times_{SO(n)}\Lambda^n\Sigma^{n+1}X
$$
which leads to an analysis of the cohomology of labelled spaces $E(\eta; X)$, paying special attention in the case when the space $X$ is a sphere as it is the main tool for calculations.

\medskip
In the case of mod-$2$ coefficients one first considers the bundle $\eta$ above and the induced fibration on bundle-configurations
\begin{equation}\label{ee}
C(S^n;S^q)\simeq \Lambda^nS^{n+q}\lra E(\eta;S^q)\lra BSO(n+1)
\end{equation}
In order to analize the fiber one considers the fibration $\Omega^nS^{n+q}\ra \Lambda^nS^{n+q}\ra S^{n+q}$ whose mod-$2$ spectral sequence collapses giving the isomorphism
$$
H^*(\Lambda^n S^{n+q};\ff_2)\cong H^*(S^{n+q};\ff_2) \otimes H^*(\Omega^n S^{n+q};\ff_2)
$$
where the second factor on the right can be described from duals of Araki-Kudo-Dyer-Lashof operations as $H^*(\Omega^n S^{n+q};\ff_2)= H^*(S^q;\ff_2) \otimes B_{n,q}$, where $B_{n,q}$ is the exterior algebra given as the kernel of 
$$
H^*(\Omega^nS^{n+q};\ff_2)\lra H^*(S^q;\ff_2)
$$
With this description, the $E_2$-term of the spectral sequence of the fibration (\ref{ee}) has the form
$$
E_2^{*,*}=H^*(BSO(n+1);\ff_2)\otimes H^*(S^{n+q};\ff_2)\otimes H^*(S^q;\ff_2) \otimes B_{n,q}
$$
converging to $H^*(E(\eta;S^q);\ff_2)$. All the differentials for this spectral sequence are described in \cite{BodigheimerCohenPeim01} and it turns out that there is an isomorphism of $H^*(BSO(n+1);\ff_2)$-modules
\begin{eqnarray*}
H^*(E(\eta;S^q);\ff_2)&\cong [H^*(BSO(n+1);\ff_2)\otimes (B_{n,q}\oplus x_qx_{q+n}B_{n,q})] \oplus\\ & [H^*(BSO(n);\ff_2)\otimes x_qB_{n,q}]
\end{eqnarray*}

Now, specializing to case of the bundle $\eta:BSO(2)\ra BSO(3)$ we get that for $k\geq 2$, there is an isomorphism of $H^*(BSO(3);\ff_2)$-modules
$$
H^*(\omcg^{2k};\ff_2)\cong H^*(BSO(3);\ff_2)\otimes H^*(F_{2k}(S^2)/\Sigma_{2k};\ff_2), 
$$
where the mod-$2$ cohomology of $F_{2k}(S^2)/\Sigma_{2k}$ is isomorphic to \\ $H^*(B_{2k};\ff_2)\oplus H^{*-2}(B_{2k-2};\ff_2)$ 
as vector spaces (see Section \ref{sub3.1} above). For a number odd of points the isomorphism has the form
$$
H^*(\omcg^{2k-1};\ff_2)\cong H^*(BSO(2);\ff_2)\otimes H^*(B_{2k-2};\ff_2). 
$$
for $2k-1\geq 3$ with anomalies for $2k-1=1$.

\medskip
Let $p$ be an odd prime and consider the sign representation $\ff_p(\pm 1)$ as a trivial $\omcg^k$-module. From the first comments on this section the homology $H_*(\omcg^k;\ff_p(\pm 1))$ is isomorphic to $\wt{H}_{*+k(2q+1)}(D_k(\eta;S^{2q+1});\ff_p)$, for $\eta$ given as above. In this new setting the Serre spectral sequence of the bundle (\ref{ee})
$$
C(S^2;S^{2q-1})\lra E(\eta;S^{2q-1})\lra BSO(3)
$$
collapses with mod-$p$ coefficients. From this and the fact that the mod-$p$ cohomology of $BSO(3)$ is that of $BS^3$ we get the isomorphism of $H^*(BS^3;\ff_p)$-modules
$$
H^*(E(\eta;S^{2q-1});\ff_p)\cong H^*(BS^3;\ff_p)\otimes H^*(\Lambda^2S^{2q+1};\ff_p)
$$
Considering the dual isomorphism in homology and noting this preserves filtration one gets the isomorphism of vector spaces 
$$
\wt{H}_*(D_k(\eta;S^{2q-1});\ff_p)\cong H_*(BS^3;\ff_p)\otimes \wt{H}_*(D_k(S^2;S^{2q-1});\ff_p).
$$
Moreover, for $k\geq 3$ there is an isomorphism
$$
H_*(\omcg^k;\ff_p(\pm 1))\cong \bigotimes_q H_q(BS^3)\otimes \wt{H}_{*-1-k(2q-1)}(D_k(S^2;S^{2q-1});\ff_p)
$$
If $p-1\geq k\geq 3$ then $H_*(\omcg^k;\ff_p(\pm 1))$ is trivial. The cases when $k=p$ and $k=p+1$ are also considered, see Corollary 9.2 in \cite{BodigheimerCohenPeim01}.

\medskip
Let $p$ be an odd prime and consider $\ff_p$ as a trivial $\omcg^k$-module. Since $\omcg^3$ is isomorphic to $\Sigma_3$, the symmetric group on $3$ letters, one has for $p>3$
$$
H^*(\omcg^3;\ff_p)\cong \begin{cases} \ff_p, & *=0\\0,&*>0            \end{cases}
$$
Now, if $p=3$ then $H^*(\omcg^3;\ff_p)\cong \Lambda[u]\otimes \ff_p[v]$, for generators $u,v$ of degrees $3$ and $4$, respectively (\cite{AdemMilgram94}). 

\medskip
As before consider the fibration $C(S^2;S^{2q})\ra E(\eta;S^{2q})\ra BSO(3)$. Since its Serre spectral sequence collapses with mod-$3$ coefficients we get an isomorphism of $H^*(BS^3;\ff_3)$-modules
$$
H^*(E(\eta;S^{2q});\ff_3)\cong H^*(BS^3;\ff_3)\otimes H^*(\Lambda^2S^{2q+2};\ff_3)
$$

Considering the dual isomorphism in homology and noting this preserves filtration one gets the isomorphism of vector spaces 
$$
\wt{H}_*(D_k(\eta;S^{2q});\ff_3)\cong H_*(BS^3;\ff_3)\otimes \wt{H}_*(D_k(S^2;S^{2q});\ff_3).
$$
Moreover, from the isomorphism $H_*(\omcg^k;\ff_3)\cong \wt{H}_{*+2qk}(D_k(\eta;S^{2q});\ff_3)$ one finally gets
$$
H_*(\omcg^k;\ff_3)\cong \bigotimes_q H_q(BS^3;\ff_3)\otimes \wt{H}_{*-1-2kq}(D_k(S^2;S^{2q});\ff_3)
$$
From this expression the Euler-Poincar\'e series for the mod-$3$ homology of $\omcg^k$ can be obtained (see \cite{BodigheimerCohenPeim01}, p.16). Moreover, these results agree with previous calculations made in \cite{BensonCohen91}.

\medskip
For $p\geq 5$ the answers are not clean as one might want. As before, the answers come from the analysis of the spectral sequence of the fibration
$$
C(S^2;S^{2q})\ra E(\eta;S^{2q})\ra BSO(3)
$$
which has a single non-zero differential. This analysis involves the device of stable Hopf invariants for showing that certain homology classes for $\Omega^2S^{2q+2}$ have non-trivial image in the homology of $E(\eta;S^{2q})$. At the end, the cohomology of $E(\eta;S^{2q})$ is given in terms of two vector spaces $A_{2q},U_{2q}$ defined from the cohomology of the cochain complex $V_{2q}$, which is isomorphic as an algebra to the product $H^*(S^{2q+2})\otimes H^*(\Omega S^{2q+1})\otimes H^*(S^{4q+1})$. See Section 8 in \cite{BodigheimerCohenPeim01} for more details.

\subsection{Applications to the genus two mapping class group}

The above methods apply to give the precise cohomology of the genus two mapping class group as follows.
There is the classic short exact sequence due to Bergau-Mennicke (\cite{BergauMennicke60}), and Birman-Hilden (\cite{BirmanHilden71})
$$
1 \to \bz/2\bz\ra \omcg_2\ra \omcg^6 \to 1.
$$
which is gotten from applying the fundamental group to a natural bundle analogous to the 
ones above. This bundle is ``twisted" by the hyperelliptic involution which is taken into account in two separate steps.

\medskip
First, let the symmetric group $\Sigma_k$ act naturally on $F_k(M)$, by the sign of a permutation on $S^1$ (thus by multiplication by $\pm 1$) and consider the orbit space $F_k(M) \times_{\Sigma_k} S^1$. In order to have an action of the unitary group $U(2)$ on $F_k(S^2) \times_{\Sigma_k} S^1$ we first consider the group $S^3 \times S^1$ acting on $F_k(S^2) \times S^1$ by the recipe
$$
(\rho, \alpha)((z_1,\cdots, z_k), \beta) = ( (\rho(z_1),\cdots, \rho(z_k)), \alpha^2\cdot \beta),
$$
where $\rho(z_i)$ denotes the projection $p:S^3 \to SO(3)$ with $SO(3)$ acting by rotations on $S^2$.
The reason for using the action of $\alpha$ as a square stems from the hyperelliptic involution. Observe that the central $\mathbb Z/2 \mathbb Z$ defined by 
$$
(-1,-1) \in S^3 \times S^1
$$
acts by fixing every point in $F_k(S^2)\times S^1$ because $-1\in S^3$ is sent to $1\in SO(3)$, and $(-1)^2=1\in S^1$. Thus the action of $S^3 \times S^1$ on $F_k(S^2) \times S^1$ descends to an action of the central product $S^3 \times_{\mathbb Z/2 \mathbb Z} S^1 = U(2)$ on $F_k(S^2) \times_{\Sigma_k} S^1$. Form the Borel construction 
$$
EU(2) \times_{U(2)}F_k(S^2) \times_{\Sigma_k} S^1.
$$
\
If $k = 6$, then this space is a $K(\pi,1)$ where $\pi = \Gamma_2$. In addition, if $k \geq 3$, the space is a $K(\pi,1)$, where $\pi$ is the centralizer of the hyperelliptic involution for all $k = 2t \geq 2$. Moreover, there is a fibration
$$
EU(2) \times_{U(2)}F_k(S^2) \times_{\Sigma_k} S^1 \lra ESO(3) \times_{SO(3)}F_k(S^2)/{\Sigma_k},
$$
induced by the natural map $U(2) \to SO(3)$ with fibre $\rp{\infty}$, realizing the exact sequence of Birman-Hilden/Bergau-Mennicke above. The cohomology of these spaces are then worked out via classical methods.

\medskip
For example, the $2$-torsion in the integral cohomology of $\Gamma_2$ has a direct summand of $\bz/8$ occurring in all strictly positive dimensions which are multiples of $4$ while the remaining $2$-torsion is all of order exactly $2$. Furthermore, the mod-$2$ cohomology is closely tied to the Brown-Gitler spectra (\cite{Cohen93}).

\medskip
A curious point about classifying spaces arises in this setting. Name\-ly, the above constructions give a map
$$
B\omcg_2 \lra BU(2)
$$
which induces an epimorphism in mod-$2$ homology. However it is not hard to check that this map is not induced by a homomorphism
$$
\omcg_2\lra U(2).
$$
Observe that $U(2)$ is the maximal compact subgroup of $Sp(4,\br)$. Thus it is possible that this fibration is induced by the natural homomorphism 
$$
\omcg_2\lra Sp(4,\br)
$$
to give a map on the level of classifying spaces $B\omcg_2\ra BSp(4,\br)$.

\medskip
Some related results have been produced in work by D. Petersen \cite{Petersen}. For example, Petersen was able to
work out the rational cohomology of the genus two mapping class group with coefficients in the exterior powers
of the defining $4$-dimensional symplectic representation. Other related interesting work is given in  \cite{Petersen3} as well as joint work with O. Tommasi \cite{Petersen2} for genus $0$ and $2$ surfaces with marked points.

\subsection{The pointed mapping class group of $T^2$}
The \res{pointed mapping class group} is the group $\omcg_g^{k,*}$ of path-components in $\topp^+(S_g;k,*)$, the group of pointed, orientation-preserving homeomorphism of $S_g$ which leaves a chosen point $*$ fixed and leaves a set $Q_k$ invariant.

\medskip
For $g=1$ and $k\geq 2$, the construction (\ref{borel}) is a $K(\omcg^k_1,1)$ space. For $T^2=S^1\times S^1$ one has an isomorphism $\omcg_1\cong SL(2,\bz)$ given by the action on $H_1(T^2)$ of any homeomorphism of $T^2$. The inclusion $T^2\subset \diff_0(T^2)$ is a homotopy equivalence (\cite{EarleEells69}) and so there is an exact sequence
$$
1\lra \diff_0(T^2) \lra \diff(T^2)\lra SL(2,\bz)\lra 1
$$
which is  split via the action
$$
SL(2,\bz)\times T^2\lra T^2, \;\; (M,(u,v))\mapsto (u^a v^b,u^c v^d)
$$
where $M=\Big(\begin{array}{cc}
a & b\\
c & d
\end{array}\Big)$, with $M\in SL(2,\bz)$. Note this action fixes the point $(1,1)$ and thus $SL(2,\bz)$ acts on $T'=T^2\baca (1,1)$. This implies that the action of the fundamental group of $BSL(2,\bz)$ on the cohomology of the fibre in
$$
B(S^1\times S^1)\lra B\diff(T^2)\lra BSL(2,\bz)
$$
is exactly the symmetric powers of the tautological representation of $SL(2,\bz)$. G. Shimura (\cite{Shimura71}) worked out the classical ring of modular forms by working out the cohomology of $SL(2,\bz)$ with coefficients in the symmetric powers of the tautological two dimensional representation over the reals, exactly the real cohomology of the bundle above. Subsequently, Furusawa, Tezuka, and Yagita (\cite{FurusawaTezukaYagita88}) worked out the cohomology of $B\diff(S^1\times S^1)$ in terms of the ring of modular forms as this last computation is exactly Shimura's computation with real coefficients.

\medskip
The method for working out the homology of the pointed mapping class group has several features which are described next. Since the action of $SL(2,\bz)$ on $T^2$ fixes $(1,1)$ the group $SL(2,\bz)$ acts on the  pointed mapping space $\map_*(T^2, S^n)$. As explained above, the cohomology of
$$
ESL(2,\bz) \times_{SL(2,\bz)} \map_*(T^2, S^n)
$$
gives the cohomology of the pointed mapping class group with marked points in case $n$ is even, and the cohomology with coefficients in the sign representation if $n$ is odd (with a few exceptions for $0$, and $1$ punctures).

\medskip
The cohomology groups of $\map_*(T^2, S^{2t+1})$, namely, $n$ odd, follows directly as was described above. In this case, the rational cohomology is 
$$
E[a]  \otimes \mathbb Q[ x_1,x_2]
$$
where $|a| = 2t-1$, and $|x_i|= 2t$ arising from the natural fibration
$$
\map_*(T^2, S^{2t+1}) \lra \Omega(S^{2t+1}) ^2
$$
with fibre $\Omega^2S^{2t+1}$. The spectral sequence of the fibration 
$$
ESL(2,\bz) \times_{SL(2,\bz)} \map_*(T^2, S^n) \lra BSL(2,\bz)
$$ 
collapses, and the resulting answer looks like a direct sum of $2$ copies of the results of Shimura's computation with degree shifts. More precisely, the vector space $M^0_{2n}$ denotes the vector space of modular cusp forms of weight $2n$ (using Shimura's weight convention) based on the
standard $SL(2,\mathbb Z)$-action on the upper $\frac{1}{2}$-plane $\mathbb H$. 

\medskip
Recall that $\br(\pm1)$ denotes the sign representation. The main result here is that if $ k \geq 2$, then
$$
H^i(\Gamma_1^{k,*}; \mathbb R(\pm1))=\begin{cases}
M^0_{2n+2} \oplus \br, &  \  \mbox{if}  \  k = 2n, \ i = 2n+1  \\
M^0_{2n+2} \oplus \br, &  \    \mbox{if}  \  k = 2n+1, \ i = 2n+1  \\
0, & \ \mbox{otherwise}.  \\
\end{cases}
$$
In this setting, the cohomology of the mapping class group with marked points in the {\it unpointed version}
$\Gamma_1^{k}$ admits a similar description in case $k \geq 2$ as follows.
$$
H^i(\Gamma_1^{k}; \mathbb R(\pm1))=\begin{cases}
(M^0_{2n} \oplus \br) \oplus(M^0_{2n} \oplus \br), & \   \mbox{if}  \  k = 2n, \ i = 2n+1\\
M^0_{2n+2} \oplus \br, &   \mbox{if}\ k = 2n+1, i = 2n+1\\
M^0_{2n+2} \oplus \br, &      \mbox{if}\ k = 2n+1, i = 2n+3\\
0 & \ \mbox{otherwise}.  \\
\end{cases}
$$

 \
On the other hand, the case of the trivial representation takes more work stemming from the fact that the cohomology of the pointed mapping space $\map_*(T^2, S^{2t})$ supports an unpleasant differential in the cohomology Serre spectral sequence from the ``fibre" to the ``middle of the spectral sequence". Keeping track of the answers gives a more intricate result depending on choices of certain partitions. However, the answers stabilize in the following sense: fix the homological dimension, and let the number of punctures increase, then the limit of these groups stabilize with the following simple answer.
\

\begin{enumerate}
\item $H^{4j}(\Gamma_1^{k,*}; \mathbb R)$ for $j \geq 0$:

            \begin{enumerate}
            \item $H^{0}(\Gamma_1^{k,*}; \br)= \br$, for all  $k \geq 0$.
                \item If   $k \geq  4j$ with $j > 1$ , then $H^{4j}(\Gamma_1^{k,*}; \mathbb R)$ is isomorphic to 
                $$ (M^0_{2j+2} \oplus \br) \oplus (M^0_{2j} \oplus  \mathbb R) \oplus (M^0_{2j} \oplus  \br).$$ 

       \end{enumerate}

\item $H^{4j+1}(\Gamma_1^{k,*}; \mathbb R)$ for $j >0$:
            \begin{enumerate}
            \item  If $k \geq 8j + 5$, for $j \geq 1$, then $H^{4j+1}(\Gamma_1^{k,*}; \mathbb R)$ is  isomorphic to 
            $$(M^0_{2j+2} \oplus \mathbb R) \oplus (M^0_{4j+4} \oplus  \mathbb R) \oplus (M^0_{4j+2} \oplus  \mathbb R)$$ 
            \end{enumerate}

\item $H^{4j+2}(\Gamma_1^{k,*}; \mathbb R)$ for $j \geq 0$:
            \begin{enumerate}
            \item  If $j \geq 0$ then this group is ${0}$.
            \end{enumerate}

\item $H^{4j+3}(\Gamma_1^{k,*}; \mathbb R)$ for $j >0$:
            \begin{enumerate}
            \item  If $k \geq 8j + 9$, for $j \geq 1$, then $H^{4j+3}(\Gamma_1^{k,*}; \mathbb R)$ is  isomorphic to 
            $$(M^0_{4j+6} \oplus \mathbb R) \oplus (M^0_{4j+4} \oplus  \mathbb R).$$             
\end{enumerate}
\end{enumerate}

The ranks of these groups are sometimes given by ranks of Jacobi forms as computed by Eichler and Zagier (\cite{EichlerZagier85}) and it is unclear whether this is an accident. In addition, there is another curious possible ``accident" described next. 

\medskip
The cohomology of the $3$-stranded braid group $B_3$ and $SL(2,\bz)$ with coefficients in  $\mathbb Z[ x_1,x_2]$ given by the symmetric powers of the tautological two dimensional
symplectic representation has been worked out. The integral results are at the interface of two seemingly different subjects (\cite{CallegaroCohenSalvetti13a},\cite{CallegaroCohenSalvetti13b}). 

\medskip
On one-hand, the torsion free summand of the integer cohomology corresponds to the ring of modular forms as calculated by Shimura \cite{Shimura71}, and described above. On the other hand, the $p$-torsion summand is given in terms of the cohomology of ``Anick fibrations" which have been used to bound the order of the $p$-torsion in the homotopy groups of spheres for $p>3$.

\subsection{Genus $g\geq 2$ surfaces}

The above methods admit a framework to surfaces for genus greater than $1$ which is explained in this section, but with far fewer results.

\medskip
Consider the projection map 
$$
E\diff(S)\times_{\diff(S)} F_k(S)/\Sigma_k \lra B\diff(S)
$$
with fibre $F_k(S)/\Sigma_k$. The Birman exact sequence
$$
1 \to \pi_1F_k(S_g)/\Sigma_k\lra \omcg_g^k \to \omcg_g \lra 1
$$
is obtained by applying the fundamental group to this fibration in case the genus of $S_g$ is at least 2 with modifications required for genus $0$ ($k \geq 3$), and $1$ ($k\geq 2$).

\medskip
As an aside, consider bundles and stable splittings as above.
The first step here is to describe the orders of natural regular representation bundles over configuration spaces 
of surfaces 
$$
F_k(S) \times_{\Sigma_k} (\br^n)^k \lra F_k(S)/{\Sigma_k}.
$$
It is known that if $S$ is closed orientable surface, then this bundle has order exactly $4$. Thus the bundle 
$$
F_k(S) \times_{\Sigma_k} (\br^{4m})^k \to F_k(S)/{\Sigma_k}
$$
is bundle equivalent to a product bundle, and the Thom space associated to this bundle is homotopy equivalent to $\Sigma^{4mk}F_k(S)/\Sigma_k$.

\medskip
In order to handle surfaces of genus $g > 1$ it is enough
to handle 
$$
E\diff^+(S_g) \times_{\diff^+(S_g)} C(S_g;X)
$$
for $X=S^n$. That is one setting where the bundles over configuration spaces are basic. So, in order to extend these methods, the first step is to identify the $\diff^+(S_g)$-action on the cohomology of $C(S_g;S^n)$.

\medskip
The mod-$2$ homology of $C(S_g,S^n)$ is known and is in this paper. If $n = 2q+1$, then the homology of $C(S_g, S^{2q+1})$ is known, but this is computing cohomology of the mapping class groups with marked points with coefficients in the sign representation. In this case, the rational homology of $C(S_g;S^{2n+1})$ is
$$
E[x_{2n+1}, x_{2n+3}] \otimes \mathbb Q[ y_1, \cdots, y_{2g}]
$$
where the degree of $x_i$ is $i$, and the degree of $y_j$ is $2n$. The $\diff^+(S_g)$-action is given polynomials in the natural $2g$-dimensional symplectic representation.

\medskip
Thus there is a basic question of the homology of $\omcg_g$ with coefficients in the symmetric powers of the tautological symplectic representation. This also occurs  in the cohomology without the sign representation. It is natural to ask
whether there is a closer connection to various modular forms for higher genus.

\medskip
The case of $C(S_g;S^{2n})$ is more complicated and we will try to describe the results. If $n=1\mbox{mod}(2)$ then there are fibrations
\begin{equation}\label{01}
\map_*(S_g,S^{2n+2})\lra\map(S_g,S^{2n+2})\lra S^{2n+2},
\end{equation}
and
\begin{equation}\label{02}
\Omega^2 S^{2n+2}\lra  \map_*(S_g,S^{2n+2})\lra \Omega (S^{2n+2})^{2g}
\end{equation}
whose Serre mod-$2$ homology spectral sequence collapses. Away from $2$ first note that an odd sphere is an $H$-space H-space and thus there are classical homotopy equivalences
$$
S^{2n+1}\times \Omega(S^{4n+3}) \ra \Omega S^{2n+2},\qquad \Omega S^{2n+1}\times \Omega^2(S^{4n+3}) \ra \Omega^2S^{2n+2}
$$

\
In characteristic zero there is a homotopy equivalence $\Omega S^{2n+1}\times S^{4n+1}\ra \Omega^2S^{2n+2}$ and in characteristic different than 2, the Serre spectral sequence for the fibration (\ref{02}) above
has a non-trivial differential from $E_{2n+2}^{0,4n+1}\ra E_{2n+2}^{4n+2,2n}$ that makes for the complications.

\medskip
In any case, the representations of the mapping class group which arise naturally are either the symmetric powers or exterior powers in the natural $2g$-dimensional representation. In this sense, the computations are similar to genus $1$.

\section{A taste of stable results}

The purpose of this section to draw a short comparison with the beautiful results of Madsen-Weiss \cite{MadsenWeiss07},
S.~Galatius \cite{Galatius04}, and  E.~ Looijenga \cite{Looijenga96} concerning similar spaces suitably stabilized.

\medskip
Consider the affine Grassmannian of oriented flat $d$-planes in $\mathbb R^n$ denoted $AG_{n,d}$.
These are standard flat planes not necessarily through the origin. Classically, form the colimit over $n$ 
with $d = 2$. Madsen-Tillmann exhibit a map $B\Gamma_{\infty} \to \Omega^{\infty}_0\mathcal A_{\infty,2}^+$ \cite{MadsenTillman01}.
Madsen-Weiss prove that this map induces an isomorphism in homology with any trivial coefficients \cite{MadsenWeiss07}. One consequence is their proof of Mumford's conjecture: stably the rational cohomology is a polynomial ring on generators $\kappa_i$ in degree $2i$.

\medskip
Namely, Madsen-Weiss worked out the homotopy type of the completion of the classifying space
for the stable mapping class group. One of their beautiful results is that this process gives
a familiar space $\Omega^{\infty}_0\mathcal A_{\infty,2}^+$.  A lucid description is given in A. Hatcher's survey \cite{Hatcher}. Galatius calculated the homology of $\Omega^{\infty}_0\mathcal A_{\infty,2}^+$ with coefficients in a finite field  \cite{Galatius04}.

\medskip
Ivanov proved (\cite{Ivanov93}) a stability theorem for the homology of the mapping class group with coefficients in symmetric powers of certain symplectic representations. Looijenga was able to use the ``hard Lefschetz theorem" in conjunction with this stable result to work out the characteristic zero cohomology of the stable mapping class group with coefficients in the symmetric powers of the tautological representation suitably stabilized. Namely, Looijenga worked out the cohomology of the stable mapping class group with coefficients in the symmetric powers of the tautological representation. 

\medskip
Analogous, but different results were obtained in work of Ebert-Randal-Wlliams \cite{Ebert.Williams}, O. Randal-Williams \cite{Randal-Williams}, and N. Kawazumi \cite{Kawazumi, Kawazumi2,Kawazumi3}.  Namely, structures arising from symmetric as well as exterior powers of the tautological representation are recurring in both the stable, and non-stable settings.

\section{The non-orientable setting}\label{sec6}
Let $N_g$ denote the non-orientable genus $g$ surface. The \res{mapping class group} of $N_g$ with marked points is the group of isotopy classes in $\topp(N_g;k)$, the group of homeomorphisms which leave a set of $k$ points in $N_g$ invariant, that is $\nmcg^k_g=\pi_0\topp(N_g;k)$. As in the orientable case there is the pure version $P\nmcg^k_g$.

\medskip
The study of non-orientable mapping class groups parallels the study of its orientable counterpart. For every non-orientable surface $N_{g,n}$ there is a double cover $p:S_{g-1,2n}\ra N_{g,n}$ which is used to obtain a set of generators for $\nmcg_g$ from that of $\omcg_g$ via symmetric homeomorphisms of $S_g$. (\cite{BirmanHilden71,BirmanChillingworth72}). There have been several developments on this issue, in particular, in \cite{ParisSzepietowski13} L. Paris and B. Szepietowski obtain a finite presentation for $\nmcg_{g,0}$ and $\nmcg_{g,1}$, for all $g$ such that $g+n>3$. In the stable setting, there has been interest in extending results to non-orientable mapping class groups such as homological stability. After Harer's stability for orientable mapping class groups, N. Wahl (\cite{Wahl08}) proved that for $g\geq 4q+3$ the group $H_q(\nmcg_{g,n})$ is independent of $g$ and $n$. Later, this stability property was extended to the case of marked points by E. Hanbury (\cite{Hanbury08}) by showing that, for any $k\geq 0$, the group $H_q(\nmcg_{g,n}^k)$ is independent of $g$ and $n$, when $g$ is sufficiently large compared to $q$.

\medskip
Recall that a group $G$ of finite virtual cohomological dimension is {\it periodic} if for some $d\neq 0$ there is an invertible element $u\in \hat{H}^d(G;\bz)$. Cup-product with this element gives a periodicity isomorphism $\hat{H}^i(G;\bz)\cong \hat{H}^{i+d}(G;\bz)$. Similarly, $G$ is $p$-periodic for a prime $p$, if the $p$-primary component $\hat{H}^*(G;\bz)$ contains an invertible element of non-zero degree $d$. Here $d$ is referred to as the period ($p$-period) of the group $G$. It is known that $\omcg_g$ is $2$-periodic for all $g>0$ and $\omcg_g$ is $p$-periodic if and only if $g$ and $p$ satisfy certain relations which show that $p$-period depends on $g$. It turns out that for marked points $\omcg^k_g$ is periodic and the period is always $2$ \cite{Lu01}. For the non-orientable mapping class group, it is known that $\nmcg_g$ is $p$-periodic whenever $p$ is odd and $g$ is not equal to $2\mod p$. See \cite{HopeTillman09}, where the double cover mentioned earlier has a special role on the results. The case with marked points is still open.

\medskip
Finally, it is worth mentioning that from the theory developed in the last sections the study of non-orientable mapping class groups with marked points is concerned with the spaces $E(\eta;k)/\Sigma_k$ and $E(\eta;X)$, for a suitably chosen bundle $\eta$, given in terms of configuration spaces of the surface $N_g$. Moreover, the theory developed so far shows that information can be obtained by the study of the labelled space $C(N_g;S^k)$, for a sphere of even dimension.

\section{Final remarks}
A recurring theme here is the natural action of the mapping class group given by the symmetric and exterior powers of the
standard symplectic representation. Both in the non-stable and stable cases, the (co)homology of the mapping class group or braid groups are both interesting and make contact with several parts of mathematics.  This recurring theme is suggestive, but at this writing the authors know very little about the problem in general.

\bigskip\bigskip\bigskip
\hfill\
{\footnotesize
\parbox{5cm}{Fred R. Cohen\\
{\it Department of Mathematics},\\
University of Rochester,\\
Rochester,\\
NY 14627\\
{\sf cohf@math.rochester.edu}}\
{\hfill}\
\parbox{5.2cm}{Miguel A. Maldonado\\
{\it Unidad Acad\'emica de Matem\'aticas},\\
Universidad Aut\'onoma de Zacatecas,\\
Zacatecas 98000,\\
M\'exico,\\
{\sf mimaldonadoa@gmail.com}}}\
\hfill

\end{document}